\documentclass{amsart}  
\usepackage{amsfonts}

\usepackage{tikz}
\usetikzlibrary{shapes}

\usepackage{endnotes}
\let\footnote=\endnote

\title{Max Dehn, Axel Thue, and the Undecidable} 

\author{Stefan M\"uller--Stach}
\email{mueller-stach@uni-mainz.de}
\address{Johannes Gutenberg--Universit\"at Mainz, 55099 Mainz, Germany}

\date{Final version, ***.}

\begin{document}

\maketitle 

\section*{Introduction}

Dehn was not the only mathematician to develop the question that came to be 
known as the word problem (see Chapters 6 and 7). In addition to Dehn's approach 
through geometric group theory, the word problem was formulated independently by 
Axel Thue for general tree structures in 1910 and for semigroups in 1914. 

In his book with Bruce Chandler, Dehn's student Wilhelm Magnus remarked that 
Dehn and Thue knew each other and mentioned the amazing parallel between their 
discoveries:
\begin{quote}
What appears to be incidental or, if one prefers, miraculous, is the fact that 
independent of Dehn and independent of topology, a contemporary mathematician 
had begun to ask questions of the type of the word problem in combinatorial 
group theory, but in an even more general and highly abstract setting. We are 
referring to the work of Thue, who may be considered as the founder of a general 
theory of semigroups. With one widely quoted exception, this work of his is 
largely forgotten nowadays. We do not know whether Dehn was influenced by Thue, 
and we have reasons to doubt it. We know that Dehn knew Thue personally, but 
only very superficially. Dehn mentioned Thue's work on occasion, observing that 
Thue's papers dealt with combinatorial problems. But he never used them, and 
indeed there is no known direct application of Thue's work to Dehn's 
group--theoretic problems.\footnote{See \cite[p. 54]{chandler-magnus}.} 
\end{quote}
A similar statement occurs in \cite[Footnote 5]{magnus3}. There Magnus also 
remarked that Dehn's wife Toni did not recall that Max Dehn had ever mentioned 
Axel Thue in her presence although Dehn had visited Norway quite a few times and 
was skilled in the Norwegian language. Moreover, there are no known personal 
relations among the students of Dehn and Thue's only student, Thoralf Skolem.

While Dehn's work had spread quickly, we do not know when mathematicians 
became aware of Thue's work on semigroups and what we now call Thue systems. We 
only know that a paper of Emil Post from 1947 mentions that Thue's paper from 
1914 was pointed out to Post by Alonzo Church. Around 1935--1955, the word 
problem became an attractive challenge for people working in the theory of 
computation (alias recursion theory) because it was, besides the 
Entscheidungsproblem, historically one of the first genuine mathematical 
problems which appeared to be potentially undecidable\footnote{A yes/no 
decision problem for an infinite set of mathematical objects is called decidable 
(alias algorithmically solvable, or recursive) if the set $S$ of G\"odel numbers of
the involved mathematical objects is a decidable subset of $\mathbb{N}$. This 
means that there is an algorithm which has output $1$ if $n \in S$ and output 
$0$ else.}. Our investigations indicate that Alonzo Church, Emil Post and other 
people at Princeton were a major driving force in bringing the word problem and 
the theory of computation together, thus placing the heritage of Dehn and Thue 
in the right historical context.

Many problems in mathematics are accessible through computation and algorithms. 
The example of the Euclidean algorithm quite prominently shows how effective 
mathematical thinking can be in inventing algorithms. Gottfried Wilhelm Leibniz 
was the first scientist who expressed in a precise way the role of a device 
(which he called calculus ratiocinator) being able to decide about the truth of 
all reasonable statements, not necessarily restricted to mathematics, by a sort 
of logical computation. Although Leibniz's thoughts remained in an abstract 
realm, he worked on the realization of an arithmetic calculating machine during 
his whole life. The Analytic Engine conceived by Charles Babbage and Ada 
Lovelace was another attempt -- albeit unsuccessful -- towards a programmable 
computer. Finally, from about 1940 on, Konrad Zuse at Berlin and John von 
Neumann at Princeton started to construct the first fully Turing complete 
computers Z3 and ENIAC. This was of course the beginning of a success story of 
incredible impact. 
 
In the early 20th century, the notion of algorithmic computability still had no 
underlying mathematically sound theory. Nevertheless people had a pragmatic idea 
what computability was supposed to mean, i.e., to reach a result in a finite 
number of computational steps. An example for this is the formulation of 
Hilbert's 10th problem:
\begin{quote} 
Given a diophantine equation with any number 
of unknown quantities and with rational integral numerical coefficients: To 
devise a process according to which it can be determined by a finite number of 
operations whether the equation is solvable in rational 
integers.\footnote{German original in \cite{hilbert1900}: Eine diophantische 
Gleichung mit irgendwelchen Unbekannten und rationalen 
ganzen Zahlenkoeffizienten sei vorgelegt: Man soll ein Verfahren angeben, nach 
welchem sich mittels einer endlichen Zahl von Operationen entscheiden 
l\"a{\ss}t, ob die Gleichung in ganzen rationalen Zahlen l\"osbar ist.}
\end{quote}
In another direction, Hilbert started his program in proof theory (Hilbert's 
program) after 1917 to obtain a solid foundation of mathematics with a method he 
called finitistic.\footnote{The finitistic approach somehow rejects the use of 
infinitely many steps. This is related but not the same as the 
intuitionistic and constructivistic approach of Brouwer, Kronecker, and Weyl which 
Hilbert disliked. G\"odel's system $T$, or equivalently Gentzen's proof for the 
consistency of arithmetic, may also be considered as finitistic in some sense.} 
Hilbert showed a lot of optimism\footnote{See Hilbert's famous words: ''Wir 
m\"ussen wissen, wir werden wissen''.} that all metamathematical questions could 
be settled within a mathematical proof theory. In 1928, together with Ackermann, 
he posed his famous Entscheidungsproblem (decision problem). It asks for 
an algorithm in the spirit of Leibniz, which decides the provability of 
statements in first order axiomatic theories. By G\"odel's completeness theorem 
for first order logic, which he proved in his dissertation from 1929, this is 
equivalent to asking for satisfiability in all possible set--theoretic 
models. 

In 1931, G\"odel's two famous incompleteness theorems were 
discovered.\footnote{G\"odel only published the first theorem, see \cite[Vol. 
I]{goedel} and \cite{plato} for the full story. Here, completeness has a 
different meaning than in G\"odel's dissertation.} The first theorem shows that 
first--order Dedekind--Peano arithmetic is incomplete in the sense that there 
are statements which are neither provable nor disprovable but true in the 
standard model. The second theorem states that the consistency of a theory at 
least as rich as Dedekind--Peano arithmetic cannot be proved as a syntactic 
formula within the theory. Hilbert's proof theoretic program had to be modified 
in the sequel. G\"odel's proof used primitive recursive functions and the 
technique of G\"odel numberings of arithmetic statements and proofs. Tarski had 
independently shown that arithmetic truth predicates exist only outside the 
realm of the theory, which implies G\"odel's incompleteness 
theorem.\footnote{This means that the set of G\"odel numbers of true statements 
in the standard model of the natural numbers is not an arithmetic set, hence not 
even recursively enumerable.} From his correspondence with John von Neumann we 
know that G\"odel was aware of this result. However, he was not able to solve 
the Entscheidungsproblem at that time since the theory of computable functions 
was only developed in full depth after 1936.

That a given mathematical problem like Hilbert's 10th problem 
or the word problem might be undecidable was probably considered unlikely by 
most people before 1931. But after G\"odel's achievements this became a more
realistic possibility. We will see, however, that Dehn and Thue already 
realized the difficulty of the word problem very clearly around 1910. 
 
A full--fledged theory of computation emerged around 1936 through the work of 
Church, G\"odel, Herbrand, Kleene, Markov, Post and Turing. Using this, the 
undecidability of the Entscheidungsproblem and of the related Halteproblem was 
shown. See \cite{davis} for the precise history of these developments. In 
addition, Church showed the undecidability of the word problem 
for finitely generated semigroups\footnote{Church's example was not finitely 
presented, as it had infinitely many relations.} in 1937. It still took many 
more years before Post and Markov gave the first proof of the undecidability of 
the word problem for finitely presented semigroups in 1947. Another five years 
passed until the word problem for finitely presented groups was shown to be 
undecidable by William Boone and Pyotr Novikov in 1952. Two decades later, the 
undecidability of Hilbert's 10th problem was shown by Yuri Matiyasevich in 1970, 
building up on work of Martin Davis, Hilary Putnam and Julia Robinson, see 
\cite{matiyasevich1970}.

In the following text, we describe the impact that both Dehn and Thue had on 
the community of recursion theory, i.e., the theory of computation, and we shed 
some light on the period between 1936 and 1955 during which many people worked 
on proving the (un)solvability of the word problem. I would like to thank Steve 
Batterson, Martin Davis, Catherine Goldstein, Jemma Lorenat, John 
McCleary, Carl--Fredrik Nyberg--Brodda, Edmund Robertson, David Rowe, Marjorie 
Senechal, Reinhard Siegmund--Schultze, J\"orn Steuding and Marcia Tucker for 
helpful remarks and assistance. 

\section{Max Dehn and the Word Problem for Groups}


In his paper \cite{dehn1911} on the word problem, Dehn used the presentation of 
a (finitely presented) group $G$ by generators and relations.\footnote{This is 
usually denoted by $G = \langle \underbrace{s_1,...,s_n}_{\rm generators} \mid 
\underbrace{r_1,...,r_m}_{\rm relations} \rangle$. A relation $r$ is given by a 
word $r$, and the notation amounts to identifying every occurrence of $r$ with 
the trivial word, i.e., setting $r=1$. Here, a word $w$ (of length $\ell$) is a 
finite combination of generators, possibly with repetition: $w=g_1^{\pm 1} 
\cdots g_\ell^{\pm 1}$. The length of a word $w$ is denoted by $|w|$. The 
inverse $w^{-1}$ of a word $w$ is obtained by inverting all $g_i$ involved and 
reversing the order, e.g., $(g_1g_2)^{-1}=g_2^{-1}g_1^{-1}$.} As Dehn remarked, 
this concept had been studied before in detail by Walther von Dyck \cite{dyck}. 
It came up even earlier in the work of William Rowan Hamilton.

In this paper, Dehn phrased the word problem, the conjugacy problem and the 
isomorphism problem for groups. The word problem asked to decide whether a 
given word $w$ in $G$ is equal to $1$. The conjugacy problem extended this 
question to determine whether two words $w,w'$ in $G$ are conjugate and if they 
are, find $u$ such that $w'=uwu^{-1}$. The conjugacy problem implies the word 
problem, since a word $w$ is equal to $1$ if and only if it is conjugate to $1$. 
Finally, the isomorphism problem aimed to determine whether two given groups 
$G$ and $G'$ are isomorphic.

In his own words, Dehn formulated the word problem, which he called 
Identit\"atsproblem, as follows: 
\begin{quote}
Let an arbitrary element of a group by given by its composition 
out of generators. One shall provide a method which decides in a finite number 
of steps whether this element is equal to the identity or not.\footnote{German 
original in \cite{dehn1911}: 1. Das Identit\"atsproblem: Irgend ein Element 
der Gruppe ist durch seine Zusammensetzung aus den Erzeugenden gegeben. Man 
soll eine Methode angeben, um mit einer endlichen Anzahl von Schritten zu 
entscheiden, ob dies Element der Identit\"at gleich ist oder nicht.}  
\end{quote}
Dehn was aware that the word problem might turn out to be difficult for a 
general group. He wrote: 
\begin{quote}
Here we have three fundamental problems whose solution is very important and probably not possible 
without a thorough study of the subject.\footnote{German original in \cite{dehn1911}: Hier sind es 
vor allem drei fundamentale Probleme, deren L\"osung sehr wichtig und wohl nicht ohne 
eindringliches Studium der Materie m\"oglich ist.} 
\end{quote}
Hence, Dehn was aware of the difficulty of the word problem, perhaps even of 
its potential unsolvability.\footnote{Magnus in \cite[p. 55]{chandler-magnus} 
cites Dehn as follows: Solving the word problem for groups may be as impossible 
as solving all mathematical problems.} In his considerations, Dehn used 
what he called the Gruppenbild (Cayley graph\footnote{Given a finitely 
presented group $G=\langle S \mid R \rangle$ with set of generators 
$S=\{s_1,\ldots,s_n\}$, the vertices are all elements of $G$, and the 
(directed) edges connect $g$ and $gs$ for every $g \in G$ and $s \in S \cup 
S^{-1}$. The edges are usually colored.}). In figure~\ref{cayleyfigure} 
this is illustrated for the dihedral group\footnote{$D_n=\langle \sigma, \tau 
\mid \sigma^n=\tau^2=1, \tau \sigma=\sigma^{-1} \tau \rangle$.} $D_5$. In 
\cite{dehn1912}, Dehn solved the word problem for infinite surface groups, i.e., 
fundamental groups of orientable closed $2$-manifolds.\footnote{These are of the 
form $G=\langle a_1,b_1,...,a_g,b_g \mid \prod_{i=1}^g a_ib_ia_i^{-1}b_i^{-1} 
\rangle$ with one relation.} 

\begin{figure}[h!]  
\begin{tikzpicture}[rotate=90,scale=1.5, thick]
\tikzstyle{vertex}=[draw,thick,circle,fill=darkgray!25,minimum size=20pt,
inner sep=0pt]

\draw (-0*360/5:1) node[vertex] (v1) {$\tau$};
\draw (-0*360/5:2) node[vertex] (v2) {$1$};
\draw (-1*360/5:1) node[vertex] (v3) {$\tau \sigma^4$};
\draw (-1*360/5:2) node[vertex] (v4) {$\sigma$};
\draw (-2*360/5:1) node[vertex] (v5) {$\tau \sigma^3$};
\draw (-2*360/5:2) node[vertex] (v6) {$\sigma^2$};
\draw (-3*360/5:1) node[vertex] (v7) {$\tau \sigma^2$};
\draw (-3*360/5:2) node[vertex] (v8) {$\sigma^3$};
\draw (-4*360/5:1) node[vertex] (v9) {$\tau \sigma$};
\draw (-4*360/5:2) node[vertex] (v10) {$\sigma^4$};

\draw[dashed] (v2) --(v1);
\draw (v3) -- (v1);
\draw (v4) -- (v2);
\draw[dashed] (v4) -- (v3);
\draw (v5) --(v3);
\draw (v6)-- (v4);
\draw[dashed] (v6)-- (v5);
\draw (v7)-- (v5);
\draw (v8)-- (v6);
\draw[dashed] (v8) -- (v7);
\draw (v9) -- (v1);
\draw (v9)-- (v7);
\draw[dashed] (v10) -- (v9);
\draw (v10) -- (v2);
\draw (v10) -- (v8);
\end{tikzpicture}
\caption{(Undirected) Cayley graph of $D_5$ with two generators 
$\sigma=(1\,2\,3\,4\,5)$ (rotation) and $\tau=(2\,5)(3\,4)$  (reflection, 
dashed). All figures in this text were prepared by the author using 
tikz.}  
\label{cayleyfigure}
\end{figure}
The idea of this proof is described by Dehn in \cite{dehn1912}.\footnote{German 
original: Zum Beweise haben wir blo{\ss} zu zeigen, da{\ss} jeder geschlossene 
Streckenzug in dem Gruppenbild, also in dem $4p$-Eckennetz, mit einem 
Netzpolygon mehr als $2p$ Seiten gemein hat oder zweimal in entgegengesetztem 
Sinne und nacheinander durchlaufene Strecken besitzt.} 
In modern language, he wrote that one needs to prove that any non--trivial 
closed loop in the Cayley graph of a surface group $G$ of genus $g \ge 2$ 
contains more than half of the defining relations, or can be freely 
reduced.\footnote{Freely reduced words have no substrings of the form 
$x^{-1}x$ or $xx^{-1}$.} In this way, Dehn provided an algorithm (called 
Dehn's algorithm today) to solve the word problem for $G$ and some other 
groups. It can be enumerated in a quite general form:

\begin{enumerate}
\item Let any freely reduced word $w=w_0$ be given. We construct a finite 
sequence $w_0,w_1, ..., w_n$ of freely reduced words by recursion such that 
$w=w_0$ and the lengths decrease $|w_0>|w_1|>\cdots >|w_i|>  \cdots$.
\item If $w_i$ is already constructed and empty, i.e., $w_i=1$, then terminate.
\item If $w_i$ contains a subword $a$ such that for some relation $r=ab$ and 
$|a|>|r|/2$, then replace $a$ by $b^{-1}$ in $w_i$ and obtain $w_{i+1}$.
\item If not, terminate at step $i$.
\end{enumerate}

There is the notion of a Dehn presentation for groups which is sufficient for 
Dehn's algorithm to work, see \cite[p. 345]{miller2014}. An example where 
Dehn's algorithm does not apply is the genus one case, i.e., the fundamental 
group of the torus, and -- more generally -- the free abelian group 
$\mathbb{Z}^n$ for $n \ge 2$ \cite[p. 345]{miller2014}. Note that the word 
problem is nevertheless easy to solve for free abelian groups of finite rank. 

There are large classes of groups beyond surface groups for genus $g \ge 2$ to 
which Dehn's algorithm can be extended. One direction where this was successful 
is the field called small cancellation theory. It deals with (finitely 
presented) groups where the relations have small overlap. We refrain from 
presenting any definitions and refer to the books \cite{lyndonschupp} and 
\cite{sims} for an account of this theory. Historically, small cancellation 
theory was mainly developed in \cite{tartakovskii}, \cite{greendlinger}, 
\cite{lyndon} and \cite{schupp}. For example, 
in \cite{greendlinger} it is proved that a group satisfying a small 
cancellation property denoted by $C'(1/6)$ has solvable word problem.

Small cancellation is not a geometric concept. A geometric class of finitely 
presented groups where Dehn's algorithm works are word--hyperbolic 
groups\footnote{Hyperbolic groups are defined as follows. Consider the Cayley 
graph of $G$ and endow it with its graph metric. Then $G$ is word--hyperbolic, 
if the resulting topological space is hyperbolic in the sense of \cite{gromov}, 
i.e., there is a constant $\delta >0$ such that any triangle is $\delta$--thin.} 
which satisfy certain metric conditions on the Cayley graph, see \cite{gromov}. 
Small cancellation groups satisfying the $C'(1/6)$--condition are examples of 
word--hyperbolic groups. It is a theorem due to Gromov and Olshanskii that for a 
general group $G$ -- in the sense that $G$ is in some way chosen randomly -- 
the Dehn algorithm solves the word problem for $G$, see 
\cite{gromov,olshanskii}. 

For other algorithms related to the word problem see \cite{knbe} and 
\cite{todd}. The Knuth--Bendix algorithm\footnote{Donald Knuth was a 
great--great--grandstudent of Thue via Thue--Skolem--{\O}re--Hall--Knuth.}
for completing term rewriting systems can be used to solve the word problem for 
the large class of {automatic groups \cite{epstein} which contains 
word--hyperbolic groups and braid groups. The Todd--Coxeter algorithm, which is 
primarily a coset enumeration method for finite index subgroups, can also be 
applied to the word problem.

The historical survey of John Stillwell \cite{stillwell} on the word 
problem contains many examples of finitely presented groups with a 
solvable word problem. In the following table we list some of them: 

\begin{figure}[h!]
\begin{tabular}{l|l} \hline 
Type of group & Reference \\ \hline
Surface groups     & \cite{dehn1912}\\ 
Trefoil knot group & \cite{dehn1914}\\ 
Subgroups of free groups (abelian or not) & \cite{nielsen}\\ 
Braid groups & \cite{artin} \\ 
One--relator groups & \cite{magnus1} \\ 
Residually finite groups\footnote{As later observed in \cite{huberdyson} as 
well as in \cite{mostowski}.} & \cite{mckinsey} \\ 
Hypo--abelian groups & \cite{engel} \\ 
Linear groups & \cite{rabin} \\ 
Knot groups & \cite{waldhausen} \\ 
Hyperbolic groups & \cite{alonso} \\ 
Automatic groups  & \cite{epstein} \\ \hline
\end{tabular}
\end{figure}

Among these people, Engel and Magnus were students of Dehn (see chapter 
6).\footnote{A list of students of Dehn is contained in \cite{magnusmoufang}.} 
Magnus proved his Freiheitssatz\footnote{The Freiheitssatz asserts that leaving 
away at least one generator appearing in the relation induces a free subgroup 
in any one--relator group $G$.} in 1930 to treat the one--relator case. 
Amazingly, the word problem for one--relator semigroups is still open, see 
\cite{fredrik}. 

Other finitely presented groups for which the word problem has been solved are 
finite groups, polycyclic groups, Coxeter groups and finitely presented simple 
groups. We refer to the textbooks \cite{lyndonschupp}, \cite{sims} for these and 
other cases. 

\section{Axel Thue and the Word Problem for Semigroups} 

Axel Thue was a number theoretist with broad interests and he was well--known 
for his work in arithmetic far beyond Norway. He held a chair position in 
applied mathematics at Oslo from 1903 on. Some of Thue's most important work in 
number theory is concerned with diophantine equations. For example, he looked at 
integer solutions of equations $f(x,y)=c$ for a homogenous polynomial $f$ with 
integer coefficients and showed that the number of those solutions is finite, 
provided certain conditions on $f$ are valid, in particular the degree of $f$ 
needs to be at least three.\footnote{The equation defines a plane curve in the 
projective plane of the same degree with equation $f(x,y)=c \cdot z^{\deg(f)}$. 
The other conditions on $f$ which we did not mention take care that this curve 
is not rational, i.e., the image of a projective line.} Such results were later 
extended by Carl Ludwig Siegel and are the basis of finiteness conjectures in 
modern arithmetic geometry (on Siegel, see chapter 5). In the same paper 
\cite[p. 232]{thue}, published in Crelle's Journal in 1909, Thue looked at 
generalizations of Liouville's result which bounds the approximation of 
irrational algebraic numbers by rational numbers from below. Thue's results 
later were strengthened by Siegel in his 1929 dissertation under Edmund Landau 
and in 1955 by Klaus Friedrich Roth who obtained an optimal estimate. Today the 
final result is known as the Thue--Siegel--Roth theorem.\footnote{The theorem 
asserts that for every algebraic number $\alpha$ and every $\varepsilon>0$ the 
inequality $\left\lvert \alpha - \frac{p}{q}\right\rvert < q^{-2-\varepsilon}$ 
has only finitely many solutions in coprime integers $p,q$.} As 
a consequence, Roth received a fields medal during the ICM at Edinburgh in 
1958. 

Thue claimed that it happened often that he discovered results which were 
previously obtained by others. For example, he wrote in a letter to Elling Holst 
from 1902 \cite[p. xxi]{thue} that he had discovered the transcendence of $e$ 
and $\pi$ independently of Hermite and Lindemann during his time as a teacher at 
the technical college in Trondheim, i.e., between 1894 and 1902. 

Among Thue's many papers are also four quite abstract papers about trees, words, 
semigroups and term rewriting which were written in German and belong to 
mathematical logic: 

\begin{itemize}
\item ''\"Uber unendliche Zeichenreihen'' \cite[p. 139--158]{thue}  from 1906.
\item ''Die L\"osung eines Spezialfalles eines generellen logischen Problems'' 
\cite[p. 273--310]{thue} from 1910.
\item ''\"Uber die gegenseitige Lage gleicher Teile gewisser Zeichenreihen''  
\cite[p. 413--477]{thue} from 1912. 
\item ''Probleme \"uber Ver\"anderungen von Zeichenreihen nach gegebenen 
Regeln'' \cite[p. 493--524]{thue} from 1914.
\end{itemize}

The papers from 1906 and 1912 present the general theory of trees and words. 
For example, in the 1906 paper Thue proves theorems which assert that there are 
infinitely long sequences consisting of three or four letters which are 
square--free, i.e., no finite length word $B$ occurs twice as $BB$ in the 
sequence. The 1906 paper continues by showing that there is an infinite 
sequence 
\[
01101001100101101001011001101001 \cdots 
\]
in two letters which is cube--free, i.e., no finite word $B$ occurs as $BBB$. 
Thue's 1912 paper elaborates on the case of two and three symbols even more and 
classifies irreducible\footnote{A sequence is irreducible if it is 
square--free, i.e., no consecutive blocks $BB$ appear.} sequences on two letters. See 
\cite{hedlund} for all this. 

It turned out that such a sequence had already been discovered before Thue by 
Eug\`ene Prouhet in 1851 (solving the Tarry--Escot problem) and later, 
independently, by others.\footnote{Notably, Marston Morse (1921), see 
\cite{morse1921}, Kurt Mahler (1929), and the chess player Max Euwe (1929).}
The sequence shows that an infinite chess game is possible without violating 
certain chess regulations \cite{morse1938,morse-hedlund}.

Thue's paper from 1910 introduced a very general philosophical (or logical) 
problem which he phrased in a metamathematical language. In modern language, he 
considered term rewriting systems for tree--like structures. The 1914 paper is 
concerned with words (Zeichenreihen) instead of binary trees. The underlying 
algorithmic problems in the case of words are known as (Semi--)Thue systems 
or as term rewriting systems \cite[p. 181]{buechi}. 

\begin{figure}[h!] 
\begin{tikzpicture}
[level distance=10mm,every label/.style={fill=gray!30,diamond},
every node/.style={fill=gray!90,circle,inner sep=1pt},
level 1/.style={sibling distance=20mm,nodes={fill=gray!60}},
level 2/.style={sibling distance=10mm,nodes={fill=gray!60}},
level 3/.style={sibling distance=5mm,nodes={fill=gray!60}}
]
\node{p}[grow=up]     
child {node{p}
child {node[right]{p}
child {node[right,label=A]{p}}
child {node[label=C]{q}}
}
child {node[right,label=F]{q}
}}
child {node{q}
child {node[right]{q}
child {node[right]{q} child {node[right,label=E]{p}} 
child {node[right,label=A]{p}}}
child {node[left]{q} child {node[above,label=D]{q}} 
child {node[left,label=A]{p} }
}
}
child {node[label=B]{p}}
};
\end{tikzpicture}
\caption{A copy of Axel Thue's figure.}
\label{thuefigure}
\end{figure}

Let us describe some more details of Thue's work. We look at finite, binary, 
rooted trees as in figure~\ref{thuefigure} (a copy of figure 3 from \cite[p. 
275]{thue}). The outer leaves correspond to variables $A$--$F$ of a certain 
type (either of type $p$ or $q$ in figure~\ref{thuefigure}). Thue explains that 
for him there is a theory of a certain logical kind behind all this (called 
Begriffe and Begriffskategorien by Thue). In the inner nodes going to the root, 
each time two values (of type $p$ resp. $q$ in figure~\ref{thuefigure}) are 
combined by a binary operation into a new value of the indicated new type. 
Hence, going all the way to the root corresponds to the computation of a tree 
automaton which computes a value of type $p$ from the given values of the 
entry variables $A$--$F$. 

\begin{figure}[h!]
\begin{tikzpicture}
[level distance=10mm,every label/.style={fill=blue!30,diamond},
every node/.style={fill=gray!30,circle,inner sep=1pt},
level 1/.style={sibling distance=20mm,nodes={fill=gray!30}},
level 2/.style={sibling distance=10mm,nodes={fill=gray!30}}
]
\node{(A+B)+C}[grow=up]     
child {node{C}
}
child {node{A+B}
child {node[right]{B}
}
child {node[]{A}}
};
\end{tikzpicture}
\hskip5cm 
\begin{tikzpicture}
[level distance=10mm,every label/.style={fill=brown!30,diamond},
every node/.style={fill=gray!30,circle,inner sep=1pt},
level 1/.style={sibling distance=20mm,nodes={fill=gray!30}},
level 2/.style={sibling distance=10mm,nodes={fill=gray!30}}
]
\node{A+(B+C)}[grow=up] 
child {node{B+C}
child {node[right]{C}
}
child {node[]{B}}
}
child {node{A}
};
\end{tikzpicture}
\caption{Two trees symbolizing $(A+B)+C$ resp. $A+(B+C)$.}
\label{treefigure}
\end{figure}

These trees may be viewed as objects representing certain algebraic or logical 
terms, as in the following example of the associativity of addition:
\[
(A+B)+C = A + (B + C).
\]
The two trees corresponding to both sides of the equation are displayed in 
figure~\ref{treefigure}. Vice versa, a binary tree corresponds to a term. In 
summary, we see that Thue had already imagined the famous correspondence between 
trees and terms. Generalizations of this occur in Post's work on canonical 
systems \cite{post0}. 

Now, term rewriting means that such trees are transformed into other trees in 
single steps by replacing (i.e., rewriting) parts according to certain rules. 
Thue thought of this term rewriting problem as an algorithmic problem about the 
relation between two given trees $A$ and $B$ in his 1910 paper:
\begin{quote}
... so we ask in other words, whether one can find trees $C_1 C_2 \ldots C_h$, 
such that $A \sim C_1 \sim C_2 \sim \ldots \sim C_h \sim B$.\footnote{German original
in \cite[p. 280]{thue}: ... so fragen wir mit anderen Worten, ob man solche B\"aume 
$C_1 C_2 \ldots C_h$ finden kann, sodass $A \sim C_1 \sim C_2 \sim \ldots 
\sim C_h \sim B$.}
\end{quote}
Thue even claimed its possible undecidability by continuing: 
\begin{quote}
A solution of this problem in the most general case may perhaps be connected 
with unsurmountable difficulties.\footnote{See \cite{stth}. German original in 
\cite[p. 280]{thue}: Eine L\"osung dieser Aufgabe im allgemeinsten Falle 
d\"urfte vielleicht mit un\"uberwindlichen Schwierigkeiten verbunden sein.}
\end{quote}

Thue's 1910 paper essentially contains the word problem without any relations 
\cite[p. 235]{buechi}. In the 1914 paper, Thue reduced this problem from binary 
to unary trees, i.e., to words or strings of letters (Zeichenreihen) and he 
introduced (Semi--)Thue systems consisting of a finite set of words $a_1, 
\ldots, a_n$ over a given countable alphabet together with a finite set of 
operations (called productions) given by pairs of words 
$(g,h)$. Any word of the form $xgy$ with possibly empty words $x,y$ may be 
replaced by the word $xhy$ for any given production $(g,h)$. A (Semi--)Thue 
system is called Thue system, if for each production $(g,h)$ also the inverse 
production $(h,g)$ is contained in $G$. By composing words via concatenation, 
Thue systems can be viewed as semigroups with a finite presentation. This 
paper resembles a remarkable point in history where the idea of a semigroup was 
born.

With this setup, the term rewriting problem becomes the word problem for 
(finitely generated) semigroups in \cite[p. 494]{thue}. Notice that one replaces 
the question $w=1$ in the word problem for groups by $w_1=w_2$ for two words 
$w_1,w_2$ in the case of semigroups. Thue describes the problem as follows: 
assuming an arbitrary choice of given words $A$ and $B$, to find a method 
through which one can always decide after a computable number of operations 
whether any two given words are equivalent with respect to $A$ and 
$B$.\footnote{German original in \cite[p. 494]{thue}: Bei beliebiger Wahl der 
gegebenen Zeichenreihen $A$ und $B$ eine Methode zu finden, durch welche man 
nach einer berechenbaren Anzahl von Operationen immer entscheiden kann, ob zwei 
beliebige gegebene Zeichenreihen in Bezug auf die Reihen $A$ und $B$ 
\"aquivalent sind oder nicht.}

The 1914 paper is cited very frequently in the literature, for example by Post 
in 1947 \cite{post3}, whereas the three other papers are mostly unknown. Richard 
B\"uchi speculates in \cite[p. 235]{buechi} that Post might have known Thue's 
papers already in 1921, when he wrote his paper \cite{post0} on canonical 
systems which may be seen as a continuation (and extension) of Thue's ideas. But 
this is not reflected in the first sentence of \cite{post3}:
\begin{quote}
Alonzo Church suggested to the writer that a certain problem of Thue 
might be proved unsolvable.\footnote{See \cite{post3}.}
\end{quote}

\section{Dehn, Thue, and the Princeton Community}

Dehn's student Wilhelm Magnus was a faculty member in Frankfurt from 1933 to 
1938. He rejected the Nazi government in public and was suspended from office 
for this reason. During the second world war he had to work in a private 
company. In 1947 he was appointed to G\"ottingen but moved to the United States 
one year later and finally became member of the Courant Institute in 1950. 
Already in the academic year 1934/35, Magnus visited Princeton. This fact alone 
implies that in the mid 30's the Princeton community was fully aware of the word 
problem. This applies in particular to the prominent topologists Solomon 
Lefschetz, James Alexander, Ralph Fox and Marston Morse (who arrived in 1935). 
The books by Lefschetz \cite{lefschetz} (1930) and Kurt Reidemeister 
\cite{reidemeister} (1932) refer to Dehn. Alexander and Fox were 
experts in knot theory at Princeton. 

We do not know much about the dissemination of the work of Axel Thue. Although 
we suspect that his work on number theory, in particular the paper from 1909, 
had been well--known to many people, his four articles on logic were probably 
not. On the other hand, with the help of Princeton librarians we found out that 
the journal\footnote{The Kristiana Videnskabs Selskabets Skrifter, 
Mathematisk--Naturvidenskabelig Klasse I, superseeded after 1924 by the journal 
of the academy Skrifter utgitt av det Norske Videnskaps--Akademi i Oslo I, 
Matematiske--Naturvidenskabelig Klasse.} in which Thue had published his papers 
was on the shelves in Princeton university between 1894 to 1960.  

Princeton University and the Institute for Advanced Studies (IAS) play a major 
role in the development of the theory of computability and in the history of the 
word problem. While the IAS was officially independent from Princeton, there 
was significant overlap among the early mathematics faculty, including Oswald 
Veblen, John von Neumann and James Alexander (see \cite{dyson}). 

Veblen was of Norwegian descent although born in the United States in 1880. As 
professor at Princeton, Veblen spent the fall of 1913 visiting Oslo, 
G\"ottingen and Berlin \cite{batterson}. Veblen and Thue were 
both participants at the 1913 Scandinavian congress of mathematics, but we do 
not know whether they met at this occasion or at any time. 
After 1932, Veblen became a leading figure in the newly founded IAS at 
Princeton. 

In both of his academic positions, Veblen supported the hiring of people in 
seemingly remote areas like mathematical logic. 
For example, the Polish immigrant Emil Post spent the year 1920--1921 at 
Princeton as a postdoc fellow. During this time he wrote his famous article on 
canonical systems~\cite{post0}. Later he spent time at Columbia, Cornell and New 
York, often interrupted by periods in which he suffered from manic attacks. 

Without doubt, Veblen had an important impact via his student Alonzo Church who 
began studying at Princeton in 1924 and finished his dissertation under Veblen
in 1927. Church was then a postdoc in G\"ottingen and Amsterdam between 1927 
and 1928. He joined the Princeton faculty in 1929 and stayed until his 
retirement in 1967. After that he continued to teach at UCLA until 1992. Church 
had an impressive roster of students. His list of students include Boone, 
Collins, Davis, Henkin, Kleene, Rabin, Rogers, Rosser, Scott, Smullyan and 
Turing who all contributed to the theory of computation, the word problem or 
related areas of logic in some essential way. 

\section{The Rise of the Undecidable}

As we already mentioned, the year 1936 was the annus mirabilis for the theory 
of computation, also called recursion theory. It saw the birth of four notions 
of computability: the $\lambda$--calculus of Church \cite{church1936}, 
the concept of a Turing machine \cite{turing1936}, another machine concept by 
Post \cite{post1}, and the notion of partial recursive function (alias 
$\mu$--recursive functions) by Kleene \cite{kleene1}, the latter building up on 
the work of Dedekind (between 1872 and 1888), Peano (1889), Skolem (1923), 
G\"odel and Herbrand (1930--1934). Surprisingly, these four definitions are 
equivalent. It is conjectured that there is no other feasible notion of 
computability beyond them. This statement is often called Church's 
thesis.\footnote{Historically more correct it should be called 
Church--Markov--Post--Turing thesis. The relevant literature in this field is 
reprinted in \cite{davis}.} 

Equipped with a notion of a (partially defined) computable function 
$f \colon \mathbb{N} \to \mathbb{N}$, one can define recursively enumerable sets 
$S \subset \mathbb{N}$ as domains, or equivalently, as images of such maps. A set 
$S \subset \mathbb{N}$ is called decidable, if $S$ and its complement are both 
recursively enumerable, i.e., the characteristic function of $S$ is computable. 
In this way, the algorithmic (un)solvability of a logical or mathematical 
problem, i.e., the computation of the characteristic function of 
the set $S$ of G\"odel numbers associated to the instances of the problem, is 
related to the (un)decidability of $S$. 

The existence of undecidable sets is the central paradigm in this theory. First 
examples in this direction were given by sets of natural numbers related to 
undecidable problems like the Entscheidungsproblem and the 
Halteproblem.\footnote{The Halteproblem has the following simple 
interpretation. If we look at computable partial functions 
$f \colon \mathbb{N} \to \mathbb{N}$, then it is possible to define a sequence 
of Turing machines $T_n$ labeled by $n \in \mathbb{N}$ such that each computable 
partial function $f$ can be computed by at least one $T_n$. Then the set $S$ of 
all $n$ such that $T_n$ halts on input $n$ is undecidable. See \cite{davis}.} 
After it had been shown that these problems were undecidable 
(i.e., algorithmically unsolvable), people were looking for more traditional 
math problems for which undecidability could be shown. It turned out that the 
word problem for groups and semigroups was a suitable candidate. Other 
undecidable sets later occured in the negative solution of Hilbert's tenth 
problem.\footnote{This Hilbertian problem asks for an algorithm to decide 
whether a polynomial system of equations over the integers has a non--trivial 
integer solution. This problem turned out to be undecidable by showing that 
every recursively enumerable set is diophantine, i.e., the projection of the 
zero set of a system of integer polynomial equations. By applying this to an 
arbitrary undecidable set $S$, one shows that the family $X_s$ of zero sets over 
every $s \in S$ has the property that one cannot decide whether $X_s$ is empty 
or not.}

In 1937, Church announced that he could prove the undecidability of the word 
problem for a particular finitely generated semigroup which is not finitely 
presented: 
\begin{quote}
By a semigroup is meant a set in which the product of any two elements is a 
unique element of the set, the multiplication being associative but not 
necessarily obeying a law of cancellation. Consider the system of combinators, 
in the sense of Rosser (Duke Mathematical Journal, vol. 1 (1935), p. 336), 
allowing as equivalence operations $r$--conversions, $p$--conversions, and also 
the operations (allowed by Curry) of replacing $BI$ by $I$ and inversely. 
This system is a semigroup, with identity element $I$, if we take as 
multiplication the operation (introduced by Curry) which is denoted by Rosser as 
$\times$. From the relations $ab=Tb \times Ta \times B \times T$ and $T(ab)=Tb 
\times Ta \times B$ it follows that every element is expressible as a product 
formed out of the four particular elements $TI, TJ, B, T$. The semigroup thus 
has a finite set of generators, although the set of generating relations must 
apparently be infinite. There is, however, an effective process of writing out 
the series of generating relations to as many terms as desired; also an 
effective means of distinguishing generating relations from others. From the 
results of the author (American Journal of Mathematics, vol. 58 (1936), pp. 
345--363), it follows that the word problem of this semigroup is 
unsolvable. (Received April 14, 1937.)\footnote{See \cite{church1937}.}
\end{quote}
Higman type embedding theorems, which were available only much later, can be 
used to show that Church's construction can be embedded into a finitely 
presented semigroup. Post proved the undecidability of the word problem for 
finitely presented semigroups (without cancellation) in 1947 \cite{post3}. He 
mentioned that Church pointed out the 1914 paper of Thue to him. The same result 
was also proved in the same year (but independently) in \cite{markov} by A. A. 
Markov jr., the son of the famous mathematician who invented Markov processes. 
We refer to the survey article of Miller \cite{miller2014} and Rotman's book 
\cite{rotman} for more details on these and the following results.

The method Post used was to associate a (Semi--)Thue system $G_T$ to any Turing 
machine $T$ \cite{post3}. In a different way, the undecidability of the 
Halteproblem can be used to show the undecidability of the word problem for 
some (Semi--)Thue system $G_T$, see \cite[\S 33]{oberschelp} and \cite[p. 
181]{buechi}. 

Turing proved the undecidability of the word problem for semigroups admitting 
cancellation in 1950 \cite{turing1950} in an attempt to obtain the full result 
for groups. However, this result did not imply the corresponding result for 
groups, since the semigroups used in the proof at least a priori cannot be 
embedded into groups. 

The word problem for groups was successfully attacked during the following 
years (see chapters 6 and 7 for further details). 
Max Dehn died in 1952 shortly before Novikov announced his proof of the 
undecidability of the word problem for groups in \cite{novikov1}. Novikov's 
published proof of this result in \cite{novikov3} uses Turing's result from 
the 1950 paper although it employs a different method. 

Church's student William Boone independently proved the undecidability of the 
word problem for (finitely presented) groups during his thesis. His final 
results were published in \cite{boone} after a long series of six papers 
\cite{boone0}. Boone used Post's semigroup approach \cite{post2} for his proof. 
It is known that Fox and G\"odel had many conversations with Boone during this 
work. John Britton independently gave a proof in 1958 and later developed 
Boone's and other methods further \cite{britton}. There is a fascinating set of 
technical results, called Britton's lemma and Novikov's principal lemma, in 
Boone's, Britton's and Novikov's proofs which turned out to be related to each 
other, see \cite[p. 355]{miller2014} and \cite[Ch. 12]{rotman}. 

In the sequel, much simpler proofs were discovered in parallel with 
developments in group theory. One of the shortest proofs uses Higman's 
embedding theorem from \cite{higman} and we describe it in the following 
section. The same method also implies that there exists a finitely presented 
group $G$ containing isomorphic copies of all finitely generated groups having 
solvable word problem. This group $G$ then does not have a solvable word 
problem. In other words, there is no uniform algorithm for all finitely 
presented groups that have a solvable word problem.\footnote{This fact is 
called the Boone--Rogers theorem.}

We remark that there are many other properties of groups which cannot be 
recognized algorithmically, e.g., the properties of being trivial, finite, 
abelian, nilpotent, solvable, free, torsion--free, residually finite, simple or 
automatic.\footnote{This is a theorem of Adian and Rabin, see \cite[p. 
366]{miller2014} for a short proof.} It is not difficult to prove the related 
undecidability of the homeomorphism problem\footnote{This problem asks for 
deciding whether two given $n$--manifolds are homeomorphic (for $n \ge 4$). This 
was treated by A. A. Markov in 1958.} for manifolds from this. 

\section{Explicit Unsolvable Examples}

As of today, many construction principles are known that yield finitely 
presented groups for which the word problem is unsolvable (i.e., undecidable). 
One particular method is quite simple and goes back to work of Higman and 
others \cite{higman}. To obtain such an example, take the finitely generated 
(and recursively presented) group 
\[
G=\langle a,b,c,d \mid a^{-e}ba^e=c^{-e}dc^e \; \forall e \in E \rangle 
\]
where $E \subset \mathbb{N}$ is a recursively enumerable, but non--recursive 
set, i.e., an undecidable set. Then use Higman's embedding theorem \cite{higman} 
to embed $G$ into a finitely presented group $G'$ with unsolvable word problem.
Explicit examples are given in \cite{borisov} and \cite{collins}. 

An simple example of Gregory S. Tseytin \cite{tseytin,collins} for a semigroup 
with unsolvable word problem -- even in the stronger sense that on a fixed word 
$w$ the decision problem $w'=w$ for any other word $w'$ is undecidable -- is 
given by
\[ 
G=\langle a,b,c,d,e \mid ac=ca, ad=da, bc=cb,bd=db,ce=eca,de=edb, 
\]
\[
cdca=cdcae,caaa=aaa,daaa=aaa \rangle 
\]
The word in question is $w=aaa$. There is an example with $2$ generators and 
$3$ relations in \cite{matiyasevich1967}. 

\theendnotes

\end{document}